\documentclass[a4paper, oneside, 11pt]{amsart}
\usepackage[left=23mm,top=25mm,right=23mm,bottom=30mm]{geometry}
\usepackage[utf8]{inputenc}
\usepackage{amsmath,amssymb,amsthm,bbm,float,graphicx,geometry,lmodern,mathtools,parskip,setspace,subcaption,diagbox}

\usepackage{mathtools}
\mathtoolsset{showonlyrefs}

\usepackage[normalem]{ulem}

\usepackage[dvipsnames]{xcolor}

\usepackage[pagebackref=true,colorlinks=true, linkcolor=MidnightBlue, citecolor=MidnightBlue]{hyperref}
\renewcommand*\backref[1]{\ifx#1\relax \else (Cited on #1) \fi}

\usepackage{mathrsfs}
\usepackage{enumerate}
\usepackage{nicefrac}
\newcommand{\e}{\mathrm e}

\newcommand\numberthis{\addtocounter{equation}{1}\tag{\theequation}}

\usepackage{tikz}
\usetikzlibrary{backgrounds}
\usetikzlibrary{patterns}
\usetikzlibrary{positioning, shapes.geometric}

\usepackage{caption}
\usepackage{subcaption}
\graphicspath{{Images/}}
\allowdisplaybreaks

\newcommand{\N}{\mathbb{N}}

\newtheorem{theorem}{Theorem}[section]
\newtheorem{lemma}[theorem]{Lemma}
\newtheorem{corollary}[theorem]{Corollary}
\newtheorem{prop}[theorem]{Proposition}
\newtheorem{conj}[theorem]{Conjecture}

\theoremstyle{definition}

\newtheorem{remark}[theorem]{Remark}

\renewcommand{\P}{\mathbb{P}}

\newcommand{\1}{\mathbbm{1}}

\newcommand{\E}{\mathbb{E}}

\newcommand{\R}{\mathbb{R}}
\renewcommand{\d}{\mathrm{d}}

\newcommand{\Z}{\mathbb{Z}}

\newcommand*{\be}{\begin{equation}}
\newcommand*{\ee}{\end{equation}}
\newcommand*{\ba}{\begin{aligned}}
\newcommand*{\ea}{\end{aligned}}

\usepackage{bbm}

\numberwithin{equation}{section}

\title{Percolation in lattice $k$-neighbor graphs}
\date{\today}

\author{Benedikt Jahnel}
\address{Institut f\"ur Mathematische Stochastik, Technische Universit\"at Braunschweig\\
Universit\"atsplatz~2,
38106 Braunschweig, Germany, and Weierstrass Institute for Applied Analysis and Stochastics\\
Mohrenstraße 39\\
10117 Berlin\\
Germany}
\email{benedikt.jahnel@tu-braunschweig.de}
\author{Jonas K\"oppl}
\address{Weierstrass Institute for Applied Analysis and Stochastics\\
Mohrenstraße 39\\
10117 Berlin\\
Germany}
\email{jonas.koeppl@wias-berlin.de}
\author{Bas Lodewijks} 
\address{Mathematical department, Universit\"at Augsburg, Universit\"atsstraße 2, 86159 Augsburg}
\email{bas.lodewijks@uni-a.de}
\author{Andr\'as T\'obi\'as}
\address{Department of Computer Science and Information Theory, Budapest University of Technology and Economics, M\H{u}egyetem rkp. 3., H-1111 Budapest, Hungary}
\email{tobias@cs.bme.hu}

\begin{document}

\maketitle

\begin{abstract}
We define a random graph obtained by connecting each point of $\Z^d$ independently and uniformly to a fixed number $1 \leq k \leq 2d$ of its nearest neighbors via a directed edge. 
We call this graph the directed $k$-neighbor graph. Two natural associated undirected graphs are the undirected and the bidirectional $k$-neighbor graph, where we connect two vertices by an undirected edge whenever there is a directed edge in the directed $k$-neighbor graph between the vertices in at least one, respectively precisely two, directions. 
For these graphs we study the question of percolation, i.e., the existence of an infinite self-avoiding path. 
Using different kinds of proof techniques for different classes of cases, we show that for $k=1$ even the undirected $k$-neighbor graph never percolates, while the directed $k$-neighbor graph percolates whenever $k \geq d+1$, $k \geq 3$ and $d \geq 5$, or $k \geq 4$ and $d=4$.
We also show that the undirected $2$-neighbor graph percolates for $d=2$, the undirected $3$-neighbor graph percolates for $d=3$, and we provide some positive and negative percolation results regarding the bidirectional graph as well.
A heuristic argument for high dimensions indicates that this class of models is a natural discrete analogue of the $k$-nearest-neighbor graphs studied in continuum percolation, and our results support this interpretation.
\end{abstract}


\smallskip


\section{Introduction}
In the last few decades, the study of $k$-nearest-neighbor type models in continuum percolation has attracted significant attention. Häggström and Meester \cite{HaggMee96} introduced the concept of the \emph{undirected} $k$-nearest-neighbor graph, where the vertex set consists of a homogeneous Poisson point process in $\mathbb{R}^d$. In this graph, two points are connected by an edge if at least one of them belongs to the $k$ nearest neighbors of the other. A \emph{cluster} refers to a maximal connected component in this graph, and the graph is said to \emph{percolate} if it contains an infinite (or unbounded) cluster. Since the percolation probability does not depend on the intensity of the Poisson point process, the only parameters affecting the percolation behavior are $k$ and $d$. 
Häggström and Meester showed that for $k = 1$, the model does not percolate in any dimension, while for large enough $d$, percolation occurs when $k = 2$. Moreover, they established that for any $d \geq 2$ there exists a $k \in \N$ such that the graph percolates. 

This initial study of continuum $k$-nearest-neighbor graphs was later extended by Balister and Bollobás~\cite{BalisterBollobas13}, who investigated three possible senses of percolation in the model where each vertex is connected  to its $k$ nearest neighbors by a \emph{directed} edge: \emph{in-percolation} (resp.~\emph{out-percolation}), which occurs by definition whenever some point of the point process exhibits an infinite incoming (resp.\ outgoing) path ending (resp.\ starting) at it, and \emph{strong percolation}, which means that there exists a strongly connected component in the graph, i.e., a component where from any point there exists a directed path to any other point.
Additionally, they introduced another undirected graph, called the \emph{bidirectional} $k$-nearest-neighbor graph, where one connects two vertices whenever they are mutually among the $k$ nearest neighbors of each other. For the two-dimensional case, they verified percolation in the undirected graph for $k \geq 11$, in the directed graph in all the three senses (in-, out- and strong percolation) for $k \geq 13$ and in the bidirectional graph for $k \geq 15$.

Recently, Jahnel and T\'obi\'as \cite{JT22} showed that in the bidirectional graph there is no percolation for $k=2$ in any dimension, even if the underlying point process is not a Poisson point process but an arbitrary deletion-tolerant and stationary point process (satisfying some basic nondegeneracy conditions). Their proof exploits the simplicity of the structure of the bidirectional graph for $k=2$, which has degrees bounded by 2. Proving absence of percolation for $k=3,4$ seems to be out of reach at the moment even in the Poisson case, and it is also not entirely clear whether these assertions hold in all dimensions or only for $d=2$.

In this article, we introduce and analyze a discrete counterpart of the continuum $k$-nearest-neighbor model, aiming to gain a deeper understanding of its underlying structure and fundamental properties. By taking a step back and considering this discrete version we hope to shed some light on its behavior in a more controlled setting. 
For the directed $k$-neighbor graph ($k$-DnG), defined on the lattice $\mathbb{Z}^d$, each vertex is connected precisely to $k$ of its $2d$ nearest nearest neighbors, using the $\ell_1$-metric. The connections are chosen independently and uniformly (without replacement) for each vertex. By connecting nearest-neighbor pairs with undirected edges if they share at least one edge (or bidirectional edges if they share two directed edges) in the $k$-DnG, we obtain the undirected ($k$-UnG) or bidirectional ($k$-BnG) $k$-neighbor graph, respectively.

At least for high dimensions $d$ this discrete model can be expected to behave similarly to its continuum counterpart.
Indeed, if we consider the continuum nearest-neighbor percolation model with $k=2$, where the intensity of the underlying Poisson point process is such that the expected number of points in a unit ball equals one, and let $Y_i$ denote the position of the $i$-th nearest neighbor for $i\in\{1,2\}$, then H\"aggström and Meester show that $|Y_i|$ converges in probability to one as $d$ tends to infinite, as well as that the conditional distribution of $Y_i$, given $|Y_i|=r$, is uniform on the sphere $\{x\in \R^d: |x|=r\}$~\cite[Lemma $3.2$]{HaggMee96}. Moreover, the volume of the intersection of two spheres with radius $r_1,r_2\in(0.9,1.1)$, whose centers are at least $0.9$ units apart, is negligible compared to the volume of either sphere~\cite[Lemma $3.3$]{HaggMee96}. As a result, for large dimensions the continuum nearest-neighbor graph has connections at distance around $1$, which are established almost independently among different pairs of vertices. As such, an approximation of this graph on $\Z^d$ should yield similar behavior when $d$ is large. 

Let us note that the $k$-BnG can be viewed as a 1-dependent Bernoulli bond percolation model, where each edge is "open" (included in the $k$-BnG) with a probability of $p = k^2/(4d^2)$. In other words, for any given edge, it is open with probability $p$, and the events of edges being open are independent when the edges are pairwise non-adjacent. Similarly, the $k$-UnG follows the same pattern with $p = (1-k/(4d))k/d$. It is worth noting that these lattice $k$-neighbor graphs exhibit negative correlation properties, setting them apart from classical models such as the random cluster model \cite{grimmett2006random}. The presence of an edge in the models under consideration here clearly decreases the likelihood of neighboring edges being present.

Throughout this paper, our main focus revolves around investigating percolation phenomena in all three variations of the model: $k$-DnG, $k$-UnG, and $k$-BnG. Specifically, we explore the presence and absence of percolation in each variant, unraveling the intricate behavior of these models.

\subsection{Organization of the paper}
The remainder of this article is organized as follows.
We will first collect the necessary notation and formulate our main results in Section \ref{sec:set}.
In Section~\ref{sec:discussion} we then discuss our findings and mention some related conjectures and open questions.
The proofs of our main results can all be found in Section~\ref{sec:proof}. 

\section{Setting and main results}\label{sec:set}

Consider the $d$-dimensional hypercubic lattice $\mathbb Z^d$. We are interested in the percolation behavior of the {\em $k$-neighbor graph} ($k$-nG) in which each vertex independently chooses uniformly $k\le 2d$ of its $2d$ nearest neighbors. We distinguish three different types of $k$-nGs:
\begin{enumerate}
\item The {\em directed $k$-nG} ($k$-DnG) in which we open a directed edge from the vertex to each of the $k$ chosen neighbors. 
\item The {\em undirected $k$-nG} ($k$-UnG) in which we open an undirected edge from the vertex to each of the $k$ chosen neighbors.
\item The {\em bidirectional $k$-nG} ($k$-BnG) in which we open an undirected edge between two vertices whenever both choose the other one as one of their $k$ neighbors.
\end{enumerate}
We are interested in the behavior of the percolation probabilities 
\begin{equation}\label{thetaD}
\theta^{\rm D}(k,d)=\mathbb P(o\rightsquigarrow\infty \text{ in the $k$-DnG on $\mathbb Z^d$}),
\end{equation}
respectively $\theta^{\rm U}(k,d)$ and $\theta^{\rm B}(k,d)$, 
where $o\rightsquigarrow\infty$ represents the event that there exists a path along open (directed) edges from the origin to infinity. Let us note that the percolation probabilities are clearly non-decreasing functions in $k$. The dependence on the dimension $d$ is more subtle and – in contrast to many other percolation models – not monotone. 

While the outdegree in the $k$-DnG is (almost-surely) equal to $k$,  the degree of a vertex in the $k$-BnG is a binomial random variable with parameters $k$ and $k/2d$, so that the expected degree is given by $k^2/2d$ which is less than $k$, unless $k=2d$. In case of the $k$-UnG, the degree is distributed according to $k+D$, where $D$ is a binomial random variable with parameters $2d-k$ and $k/2d$, and hence, the expected degree is given by $k(4d-k)/(2d)$.

\subsection{Results for the directed $k$-neighbor graph}
As a first step, we study the cases where the occurrence of percolation is deterministic. On the one hand, choosing only one neighbor is never enough for percolation, and on the other hand, choosing a sufficiently large number of neighbors leads to the almost-sure existence of a path connecting the origin to infinity. 

\begin{prop}\label{prop:deterministic-percolation-dng}
    \begin{enumerate}[(i)]
        \item For all $d \geq 1$ it holds that $\theta^{\rm D}(1,d) = 0$.
        \item Whenever $k \geq d+1$ we have $\theta^{\rm D}(k,d) = 1$. 
    \end{enumerate}
\end{prop}

The proof can be found in Section \ref{sec:proofs-dng}. Let us now turn our attention towards the intermediate supercritical percolation phase. Here, our tools for the proof are restricted to sufficiently high dimensions $d$.

\begin{theorem}\label{thm-CoxDurrett}
    If $d \geq 4$ and $k \geq 4$ or if $d \geq 5$ and $k = 3$ we have $\theta^{\rm D}(k,d)>0$.  
\end{theorem}
The proof can be found in  Section \ref{sec:proofs-dng} and mainly relies on a variation of the technique established in \cite{CoxDurrett} plus precise estimates on the probability that two independent simple random walks meet each other and then take a common step. In fact the proof of Theorem~\ref{thm-CoxDurrett} establishes positive probabilities for the existence of infinite, self-avoiding and {\em monotone} paths, that is paths with strictly increasing $\ell_1$-norm.

Last but not least, we come to the question of monotonicity with respect to the parameters. It is also easy to see that for $d \geq k$ we have $\theta^{\rm D}(k,d) <1$ and, by Proposition \ref{prop:deterministic-percolation-dng}, $\theta^{\rm D}(k,d)=1$ for $d\le k-1$, so $d\mapsto \theta^{\rm D}(k,d)$ can be decreasing. However, in view of similar results in the continuum, see~\cite{HaggMee96}, it seems reasonable to assume that $d\mapsto \theta^{\rm D}(k,d)$ is increasing for large $d$. At least, we can deduce the following monotonicity along diagonals of the parameter space. 

\begin{theorem}\label{theorem-3322}
For all $k,d\ge1$, we have that $\theta^{\rm D}(k+1,d+1)\ge \theta^{\rm D}(k,d)$.
\end{theorem}

The proof of Theorem~\ref{theorem-3322} can be found in Section \ref{sec:proofs-dng} and is based on a coupling argument. We will see that it is essential that \emph{both} $k$ and $d$ increase by one to make the coupling work. For example, issues arise when trying to deduce a similar coupling between the settings $k=d=2$ and $k=2,d=3$. Though we expect percolation occurs in both settings, a similar coupling would give rise to a two-dimensional model where some vertices have only out-degree one. In such a case, it seems that the additional edges that appear in the $2$-DnG in two dimensions are pivotal in the sense that they enable percolation to occur whereas it may not occur without the presence of these edges. 
Still, we expect that if we restrict ourselves to parameters $k \leq d$ a similar relation between the $d$-dimensional and $(d+1)$-dimensional settings should exist. We state this in Conjecture~\ref{conj:dimensioncouple} below.

\subsection{Results for the undirected $k$-neighbor graph}

To state our results and proofs for the $k$-UnG we first introduce the following notation. 

It is clear that $\theta^{\rm U}(k,d) \geq \theta^{\rm D}(k,d)$ for any $k,d\in \N$ because any directed edge in the $k$-DnG corresponds to an undirected edge in the $k$-UnG (and at most two directed edges can correspond to the same undirected edge). This way, Theorem~\ref{thm-CoxDurrett} and Part~(ii) of Proposition~\ref{prop:deterministic-percolation-dng} also hold with $\theta^{\rm D}$ replaced by $\theta^{\rm U}$ everywhere. Moreover, we can now also deal with lower dimensions. 

\begin{theorem}\label{theorem-U22}
  We have $\theta^{\rm U}(2,2)>0$ and $\theta^{\rm U}(3,3)>0$.  
\end{theorem}

However, even in the undirected sense, $k=1$ is still not sufficient for percolation in any dimension $d \in \N$. 
\begin{prop}\label{prop-U1}
    We have that $\theta^{\rm U}(1,d)=0$ for all $d \geq 1$.
\end{prop}

\subsection{Results for the bidirectional $k$-neighbor graph}\label{sec:B}
Clearly we have 
\begin{align*}
    \theta^{\rm B}(k,d) \leq \theta^{\rm D}(k,d), 
\end{align*}
so we already know that, if each vertex chooses a single neighbor, the bidirectional model will not percolate. A stronger result holds for this model, however. The following lemma presents an upper bound for $k$ in terms of $d$ for which we can verify that the bidirectional model does not percolate. To state this lemma, let us denote by $c(d)$ the \textit{connective constant} of $\Z^d$, see e.g., \cite{grimmett2006random}, which is defined by 
\begin{align*}
    c(d) := \lim_{n \to \infty}c_n(d)^{1/n},
\end{align*}
where $c_n(d)$ is the number of self-avoiding paths of length $n$ in the $d$-dimensional hypercubic lattice that start at the origin. Via a quick subadditivity argument one can show that the limit $c(d)$ actually exists and satisfies $d < c(d) < 2d -1$. 

\begin{lemma}\label{lemma-Bnoperc}
For all $k$ such that $k(k-1)<2d(2d-1)/c(d)$ we have that $\theta^{\rm B}(k,d)=0$.
\end{lemma}
Let us note that this result can be interpreted in two ways. First, for given $d$, it guarantees absence of percolation for all sufficiently small $k$, and for example $\theta^{\rm B}(2,d)=0$ for any $d\geq 2$ (while it is clear that $\theta^{\rm B}(2,1)=1$), $\theta^{\rm B}(3,d)=0$ for any $d \geq 3$ and $\theta^{\rm B}(4,d)=0$ for any $d \geq 6$, as can be seen from the lower bounds on $c(d)$ by \cite[Table 1]{HSS} (e.g.\ $c(6)\geq 10.874038$). However, also for fixed $k$, since $2d(2d-1)/c(d)>2d$, for sufficiently large $d$, there is no percolation thanks to Lemma~\ref{lemma-Bnoperc}. This is due to the fact that, in high dimensions, it is unlikely for two neighboring vertices to pick the same connecting edge. This behavior is rather different from the case of the $k$-DnG and the $k$-UnG, where there is percolation for all $k \geq 3$ in all sufficiently high dimensions.

The proof strategy of Theorem~\ref{theorem-U22} is to verify that percolation occurs when restricting ourselves to a two-dimensional plane, see Section~\ref{sec:kungproof}. This approach is also applicable in the bidirectional model, and yields the following result.

\begin{prop}\label{prop-Bperc}
We have $\theta^{\rm B}(k,d) > 0$ whenever 
\be 
 k>d\sqrt{4\Big(1-1/c(2)\Big)}.
\ee 
\end{prop}

An application of the upper bound $c(2)\leq 2.679192495$ from~\cite{PonTitt00} for the connective constant of $\Z^2$ immediately yields the following corollary. 
\begin{corollary}\label{cor-Bperc}
We have $\theta^{\rm B}(k,d) > 0$ whenever 
\be 
 k>d\sqrt{4\Big(1-1/2.679192495\Big)}\approx 1.583355
d.
\ee 
\end{corollary}

This corollary allows us to verify percolation, e.g., for the $(2d-1)$-BnG for $d \geq 3$, for the $(2d-2)$-BnG for $d \geq 5$, for the $(2d-3)$-BnG for $d \geq 8$, and for the $(2d-4)$-BnG for $d \geq 10$. (Thus, the smallest-dimensional positive percolation result that we obtain from the corollary is that $\theta^{\rm B}(5,3) >0$.)

Note that, compared to $d$, we always need rather large $k$ to percolate. In high dimensions we can improve the ratio slightly by using that the $k$-BnG model (just as the $k$-UnG model) in any dimension $d$ is an 1-dependent Bernoulli bond percolation model where each edge is open with probability $p=k^2/(4d^2)$. By using the results from \cite{improvedLB} we obtain the following improved asymptotic result. 

\begin{prop}\label{prop-1indep}
For any $\alpha > 2 \sqrt{0.5847} \approx 1.5293$, we have $\theta^{\rm B}(\lfloor \alpha d \rfloor, d)>0$ for all $d$ sufficiently large.
\end{prop}

Of course, the $k$-UnG model is also an 1-dependent bond-percolation model, but the same approach unfortunately does not yield any new results in this case. Indeed, here each edge is open with probability $p=k(4d-k)/(4d^2)$. Thus, for $k=d$ it always holds that $p=3/4<0.8457$, while for $k>d$ (and for $d \geq 2$ even for $k \geq 4$) we already know that $\theta^{\rm D}(k,d)=1$. 

\subsection{Summary}
To close this section, let us summarise our results in Table~\ref{table-UD}. 
\begin{table}[h]
\centering
\scalebox{0.85}{
\begin{tabular}{|l|c|c|c|c|c|c|c|c|} 
 \hline
\diagbox{$d$}{$k$} & 1 & 2 & 3 & 4 & 5 & 6 & $\geq 7$ \\ \hline 
1 & no/no/no & yes/yes/yes & - & - & - & - & - \\ \hline
2 & no/no/no & no/open/yes & open/yes/yes & yes/yes/yes & - & - & - \\ \hline
3 & no/no/no & no/open/open & no/open/yes & open/yes/yes & yes/yes/yes & yes/yes/yes & - \\ \hline
4 & no/no/no & no/open/open & no/open/open & open/yes/yes & open/yes/yes & open/yes/yes & yes/yes/yes \\ \hline
$\geq 5$ & no/no/no & no/open/open & no/yes/yes & open/yes/yes & open/yes/yes & open/yes/yes & open/yes/yes \\  \hline
\end{tabular}
}
\caption{Percolation results and open cases for the BnG/DnG/UnG. For given $k,d$, `yes' means that the given graph percolates, `no' means that it does not, while `open' means that the given case is (at least partially) open.}
\label{table-UD}
\end{table}

\section{Outlook and open problems}\label{sec:discussion}
Although we managed to prove the occurrence or absence of percolation in a wide range of cases for the $k$-$\square$nG model ($\square\in\{\mathrm{D,U,B}\}$), there are still many cases where the techniques used did not provide any conclusive results, especially in low dimensions. In this section we want to comment on a few of those that we deem to be interesting and also discuss some further conjectures on the behavior of discrete $k$-neighbor graphs. 
\subsection{Directed $k$-neighbor graphs}
\begin{figure}[h]
     \centering
     \begin{subfigure}[b]{0.3\textwidth}
         \centering
         \includegraphics[width=\textwidth]{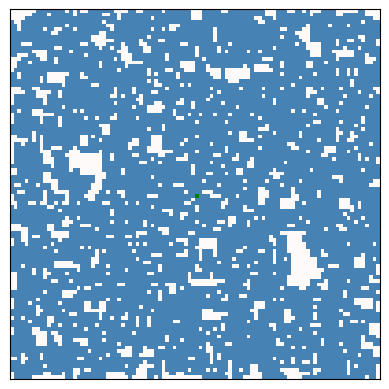}
         \caption{$n=100$}
         \label{fig:2-DnG-100x100_1}
     \end{subfigure}
     \hfill
     \begin{subfigure}[b]{0.3\textwidth}
         \centering
         \includegraphics[width=\textwidth]{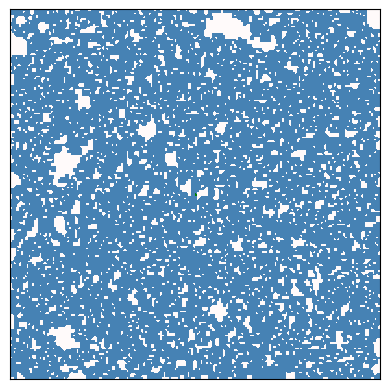}
         \caption{$n=200$}
         \label{fig:2-DnG-200x200_1}
     \end{subfigure}
     \hfill
     \begin{subfigure}[b]{0.3\textwidth}
         \centering
         \includegraphics[width=\textwidth]{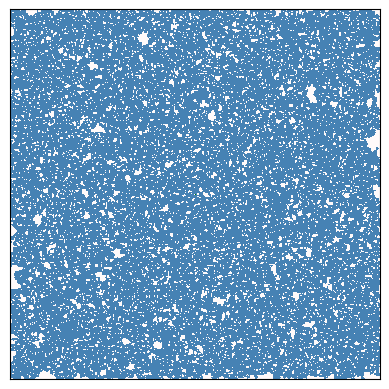}
         \caption{$n=500$}
         \label{fig:2-DnG-500x500_1}
     \end{subfigure}
        \caption{Samples of the $2$-DnG in boxes with side-length $n$. The colored vertices are the ones contained in the connected component of the origin.}
        \label{fig:2-DnG-samples}
\end{figure}

For the $k$-DnG the most pressing open question is if it percolates for $k=2$ in $d=2$ and then if we consequentially also percolate in any dimension $d$. 
According to simulations this seems to be the case, see Figure~\ref{fig:2-DnG-samples}, and, moreover, the proportion of vertices we can reach from the origin is actually quite high. This strong numerical evidence and the heuristic that it is easier for the DnG to percolate in higher dimensions suggest the following conjecture. 

\begin{conj}\label{conj:2-dng}
    The $2$-DnG percolates in all dimensions $d \in \N$. 
\end{conj}

Note that the complement of the $k$-DnG in $d$ dimensions with respect to $\Z^d$, i.e., the graph with vertex set $\Z^d$ and the edge set formed by all nearest-neighbor edges of $\Z^d$ not contained in the $k$-DnG, is in distribution a $(2d-k)$-DnG. So if Conjecture \ref{conj:2-dng} holds, then this would be an example of phase coexistence in a self-complementary (directed) graph in two dimensions. 

With regards to the monotonicity statement Theorem~\ref{theorem-3322} we expect a corresponding statement also to hold along the diagonal of the $(k,d)$-plane grid.
\begin{conj}\label{conj:dimensioncouple}
    For all $k,d\ge 2$ such that $k\leq d$, we have that $\theta^{\rm D}(k,d+1)\geq \theta^{\rm D}(k,d)$. 
\end{conj}

One possible way to prove Conjecture~\ref{conj:2-dng} would be then to verify percolation in $d=2$ and to prove Conjecture \ref{conj:dimensioncouple}, but at least for higher dimensions $d$ there might be simpler arguments making use of intersection properties of random walks to show that $\theta^{\rm D}(2,d) > 0$. Of course this would also imply that we have $\theta^{\rm U}(k,d)>0$ for all $k,d\geq2$.  

\subsection{Bidirectional $k$-neighbor graphs}

\begin{figure}[h]
     \centering
     \begin{subfigure}[b]{0.3\textwidth}
         \centering
         \includegraphics[width=\textwidth]{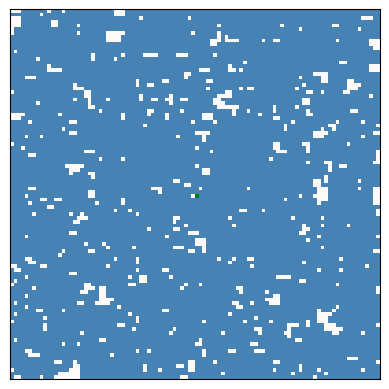}
         \caption{$n=100$}
         \label{fig:3-BnG-100x100}
     \end{subfigure}
     \hfill
     \begin{subfigure}[b]{0.3\textwidth}
         \centering
         \includegraphics[width=\textwidth]{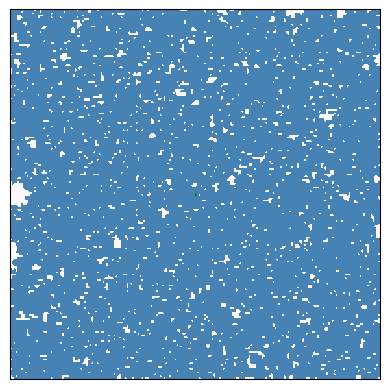}
         \caption{$n=200$}
         \label{fig:3-BnG-200x200}
     \end{subfigure}
     \hfill
     \begin{subfigure}[b]{0.3\textwidth}
         \centering
         \includegraphics[width=\textwidth]{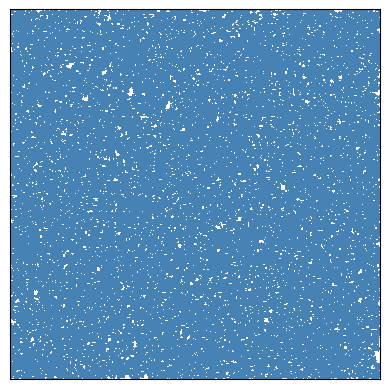}
         \caption{$n=500$}
         \label{fig:3-BnG-500x500}
     \end{subfigure}
        \caption{Samples of the $3$-BnG in boxes with side-length $n$. The colored vertices are the ones contained in the connected component of the origin.}
        \label{fig:3-BnG-samples}
\end{figure}

\begin{figure}
     \centering
     \begin{subfigure}[b]{0.475\textwidth}
         \centering
         \includegraphics[width=\textwidth]{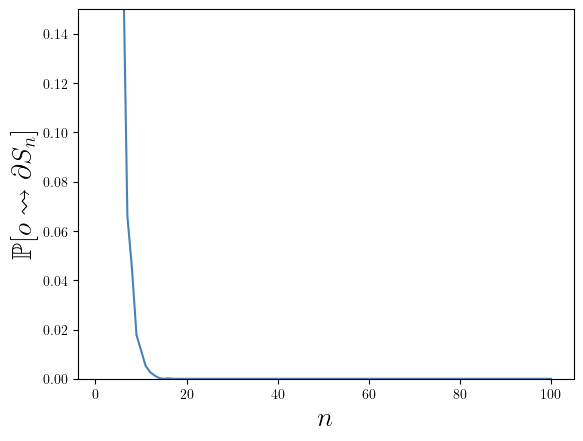}
         \caption{$d=k=2$}
         \label{fig:2-DnG-100x100}
     \end{subfigure}
     \hfill
     \begin{subfigure}[b]{0.475\textwidth}
         \centering
         \includegraphics[width=\textwidth]{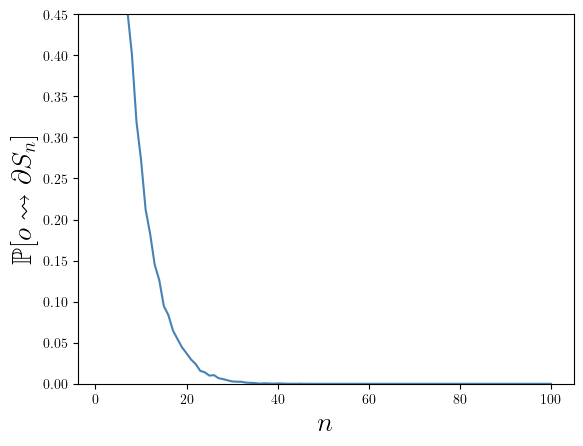}
         \caption{$d=k=3$}
         \label{fig:2-DnG-200x200}
     \end{subfigure}
     \vskip\baselineskip
     \begin{subfigure}[b]{0.475\textwidth}
         \centering
         \includegraphics[width=\textwidth]{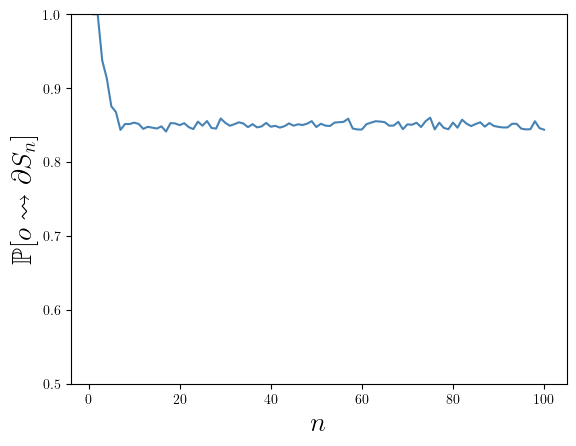}
         \caption{$d=k=4$}
         \label{fig:2-DnG-500x500}
     \end{subfigure}
     \hfill 
     \begin{subfigure}[b]{0.475\textwidth}
         \centering
         \includegraphics[width=\textwidth]{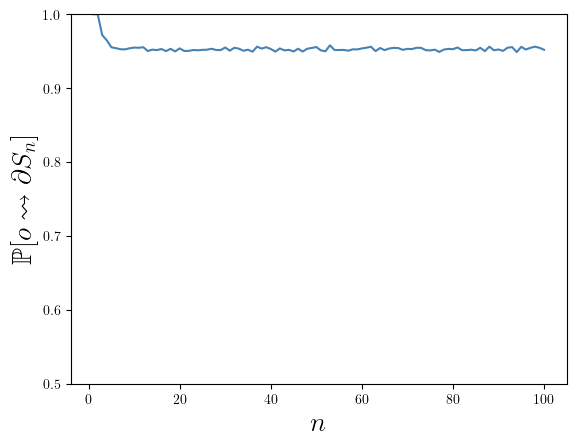}
         \caption{$d=k=5$}
     \end{subfigure}
        \caption{Probability to reach the boundary of boxes $S_n$ of varying sidelength $n$ in the $d$-BnG for $d \in \{2,3,4,5\}$ (Sample size: 10000).}
        \label{fig:BnG-probabilities}
\end{figure}

For the bidirectional model the most intriguing question seems to be what the smallest $k=k(d)$ is such that $\theta^{\rm B}(k(d),d)>0$ holds for all $d$ (or at least for all $d$ sufficiently large)? Heuristically, one would expect that we do not percolate when the expected degree of a vertex is less than $2$, because the chance of backtracking, i.e., not getting a new edge is then too large. However, as soon as we have an expected degree of at least $2$, one could imagine that this is sufficient to percolate since we usually get at least one new edge in each step while walking along a path. At least in low dimensions (which are still accessible for numerical computations) the numerical tests seem to support the following conjecture, see Figures~\ref{fig:3-BnG-samples} and~\ref{fig:BnG-probabilities}. 

\begin{conj}
    The $k$-BnG percolates in dimension $d$ if and only if $k \geq 2 \sqrt{d}$.
\end{conj}

In particular this would imply that the $d$-BnG does \textit{not} percolate for dimensions $d=2,3$ but percolates for $d \geq 4$, and that the $3$-BnG percolates in $d=2$.  

\subsection{XOR percolation}
Let us briefly consider another variant of lattice $k$-neighbor graphs, namely the {\em exclusively unidirectional $k$-nG} ($k$-XnG), in which we open an undirected edge between two vertices whenever exactly one of them chooses the other one as one of its $k$ chosen neighbors. The letter X refers to the ``XOR'' (exclusive ``or'') in the edge-drawing rule.

Although the parameter $k$ of the $k$-XnG can range between $1$ and $2d$ (just as for the $k$-BnG, $k$-DnG and $k$-UnG), actually it suffices to consider $1 \leq k \leq d$, thanks to the following lemma.
\begin{lemma}\label{lem:Xk2d-k}
For any $1 \leq k \leq 2d-1$, the $k$-XnG equals the $(2d-k)$-XnG in distribution.
\end{lemma}
\begin{proof}
The statement follows from the fact that the $k$-XnG contains precisely the undirected nearest-neighbor edges of $\Z^d$ that are included in the $k$-DnG in one direction and in the complementary $(2d-k)$-DnG in the other direction, and the same holds for the $(2d-k)$-DnG. 
\end{proof}
Of course, Lemma~\ref{lem:Xk2d-k} also holds for $k=2d$ if we define the $0$-XnG as $\Z^d$ with no edges, from which it is clear that $\theta^{\rm X}(2d,d)=0$ for all $d$, where we write $\theta^{\rm X}$ for the percolation probability of the XnG analogously to~\eqref{thetaD}. Hence, we can limit our analysis to the cases $1\le k\le d$.

\begin{figure}[h]
     \centering
     \begin{subfigure}[b]{0.3\textwidth}
         \centering
         \includegraphics[width=\textwidth]{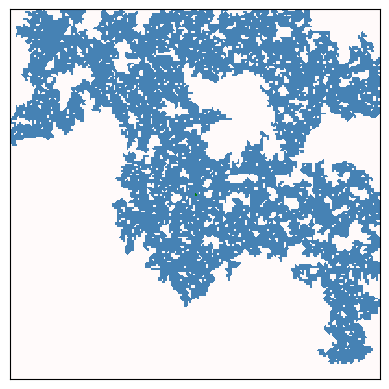}
         \caption{$n=200$}
     \end{subfigure}
     \hfill
     \begin{subfigure}[b]{0.3\textwidth}
         \centering
         \includegraphics[width=\textwidth]{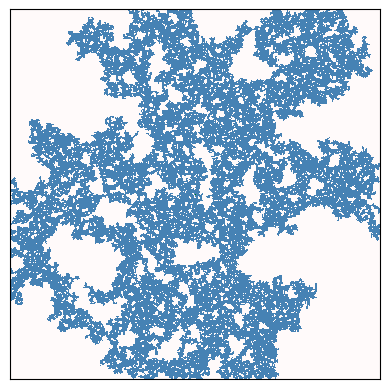}
         \caption{$n=500$}
     \end{subfigure}
     \hfill
     \begin{subfigure}[b]{0.3\textwidth}
         \centering
         \includegraphics[width=\textwidth]{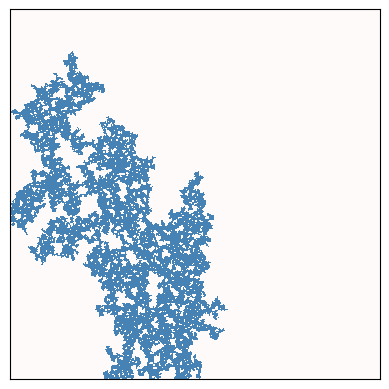}
         \caption{$n=1000$}
     \end{subfigure}
        \caption{Samples of the $2$-XnG in boxes with side-length $n$. The colored vertices are the ones contained in the connected component of the origin.}
        \label{fig:2-XnG-samples}
\end{figure}

It is also easy to see that, since the $k$-XnG is a subgraph of the $k$-UnG, we have $\theta^{\rm X}(1,d)=0$ for all $d \geq 1$, and thus by Lemma~\ref{lem:Xk2d-k}, $\theta^{\rm X}(2d-1,d)=0$ for all $d$. Moreover, the $2$-XnG is also not percolating in one dimension, so the first non-trivial case is $k=d=2$, see Figure~\ref{fig:2-XnG-samples} for some simulations. 

Since the probability of an edge being open 
in the $k$-XnG is given by $k(2d-k)/(2d^2)$, which is maximized at $k=d$ with maximum value $1/2$, the $d$-XnG seems to be critical in dimension $d$. Moreover, for $k\ge 2$, and $e_1, e_2$ edges that share a common vertex, we have that 
\be\ba 
\P(e_1,e_2\text{ open})=\frac{(2d-k)(k(4k-1)(2d-k)-k^2)}{8d^3(2d-1)}\leq \frac{k^2(2d-k)^2}{4d^4}=\P(e_1\text{ open})^2,
\ea \ee
with equality if and only if $k=d$. Thus the model features again negative correlations that decrease in $k\in\{2,\dots,d\}$ and in particular, for $k=d$, we even have independence. This, together with the fact that, for fixed $d$, the probability of an edge being open in the $k$-XnG increases, allows us to formulate the following conjecture.
\begin{conj}
For fixed $d$, $k \mapsto \theta^{\rm X}(k,d)$ is strictly monotone increasing in $k\in\{2,\dots,d\}$.
\end{conj}
Despite the above heuristic justification, it is not clear to us how to verify this conjecture. In particular, we are not aware of a coupling between the $k$-XnG and the $l$-XnG for $1 \leq k < l \leq d$, while the $k$-DnG (respectively $k$-UnG, $k$-BnG) is a subgraph of the $l$-DnG (resp.\ $l$-UnG, $l$-BnG) whenever $1 \leq k \leq l \leq 2d$.

Let us finally mention that, in view of the first-moment method for the proof of existence of subcritical regimes, see for example Lemma~\ref{lemma-Bnoperc}, we would need to establish $\P(e_1\text{ open}|e_2\text{ open})c(d)<1$, which is unfortunately not true for any $2\le k\le d$. Also none of the methods which we have used to prove existence of supercritical percolation regimes seem to apply to the $k$-XnG, so determining whether it actually percolates for certain choices of $d$ and $k$ remains a goal for our future research. 

\begin{figure}[ht]
    \centering
    \begin{subfigure}[b]{0.3\textwidth}
        \centering 
        \includegraphics[width=\textwidth]{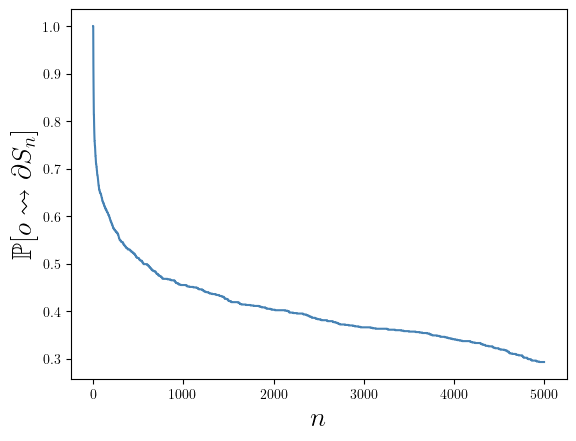}
        \caption{Probability of reaching the boundary of boxes with sidelength $n$ in the $2$-XnG in $d=2$ (Sample size=1000).}
        \label{fig:2-XnG-percolation-prob}
    \end{subfigure}
    \hfill 
    \begin{subfigure}[b]{0.311\textwidth}
        \centering 
        \includegraphics[width=\textwidth]{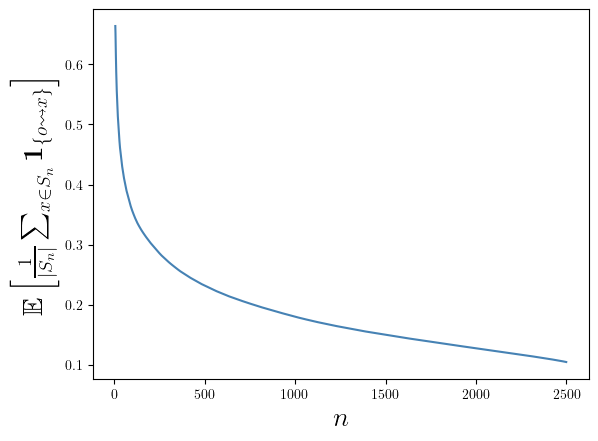}
        \caption{Proportion of vertices reached from the origin in boxes with sidelength $n$ in the $2$-XnG in $d=2$ (Sample size=1000).}
        \label{fig:2-XnG-percolation-prop}
    \end{subfigure}
    \hfill
    \begin{subfigure}[b]{0.3\textwidth}
        \centering 
        \includegraphics[width=\textwidth]{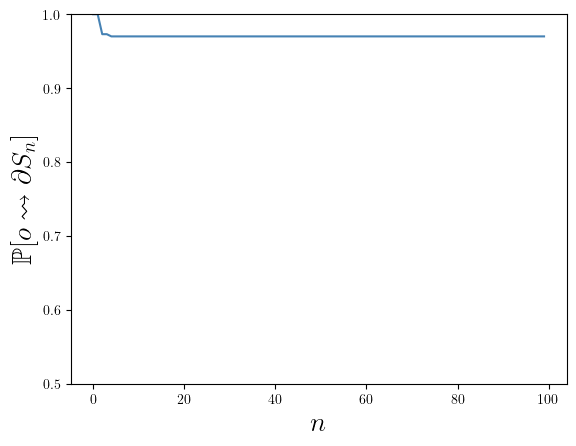}
        \caption{Probability of reaching the boundary of boxes with sidelength $n$ in the $2$-XnG in $d=3$ (Sample size=1000).}
        \label{fig:2-XnG-in-3d}
    \end{subfigure}
    \caption{Numerical experiments for the XnG in dimensions $d=2,3$.}
    \label{fig:xng-numerical-experiments.}
\end{figure}

Focussing on the case $k=2$, we have performed simulations that suggest absence of percolation for $d=2$, see Figures~\ref{fig:2-XnG-percolation-prob} and~\ref{fig:2-XnG-percolation-prop}. 

Based on this and Figure \ref{fig:2-XnG-in-3d} we at least formulate the following conjecture.
\begin{conj}
    The $2$-XnG percolates for all $d\geq 3$ but does not percolate for $d=2$. 
\end{conj}

\subsection{In-percolation and strong percolation}\label{sec-instrong}
As we already mentioned in the introduction, our notion of directed percolation~\eqref{thetaD} corresponds to out-percolation in~\cite{BalisterBollobas13}. 
In this paper, we did not develop any specific proof techniques for strong percolation or in-percolation. Nevertheless, it is clear that, for fixed $d$ and $k$, percolation in the $k$-BnG implies strong percolation, while strong percolation implies both in- and out-percolation in the $k$-DnG. Moreover, in- or out-percolation in the $k$-DnG implies percolation in the $k$-UnG, analogously to the continuum case (where these implications were mentioned in~\cite{BalisterBollobas13}). Hence, our positive percolation results for the $k$-BnG yield ones for strong percolation, and our negative (out-)percolation results for the $k$-DnG imply the absence of strong percolation. Furthermore, writing
\[ \theta^{\rm D'}(k,d) = \P(\exists \text{ an infinite directed path in the $k$-DnG ending at }o) \] 
for the analogue of $\theta^{\rm D}(k,d)$ for in-percolation, we have $\theta^{\rm D'}(1,d)=0$ since $\theta^{\rm U}(1,d)=0$ holds thanks to Proposition~\ref{prop-U1}.
\begin{conj}\label{conj:in}
    For all $d \in \N$  and $1 \leq k \leq 2d$, $\theta^{\rm D'}(k,d)>0$ if and only if $\theta^{\rm D}(k,d)>0$.
\end{conj}
\begin{remark}
Let us note that the same assertion is not true with ``with positive probability'' replaced by ``almost surely'' everywhere. Indeed, let $d+1 \leq k<2d$, then we know from Part $(ii)$ of~Proposition~\ref{prop:deterministic-percolation-dng} that $\theta^{\rm D}(k,d) =1$. However, $\theta^{\rm D'}(k,d)<1$ because the probability of the event that all the $2d$ nearest neighbors of $o$ choose only neighbors different from $o$ as endpoints of their outgoing edges is positive. 
\end{remark}

One reason why one could be led to believe that Conjecture~\ref{conj:in} is true is the following simple argument, which is based on the mass-transport principle, see e.g.~\cite[Chapter 8]{lyons_peres_2017}. Consider the function
\begin{align*}
    F(x,y;\omega) := \mathbf{1}_{\{x \to y\}}(\omega), \quad x,y\in \Z^d,
\end{align*}
where $\omega$ is a realization of the directed $k$-DnG on $\Z^d$. This function $F$ is clearly diagonally invariant with respect to lattice shifts, i.e.,
\begin{align*}
    F(x+z,y+z;\omega+z) = F(x,y;\omega).
\end{align*}
So, if we let $\overrightarrow{C_0}$ denote the \textit{out-component} of the origin and $\overleftarrow{C_0}$ denote the \textit{in-component}, then the classical mass-transport principle~\cite[p.~274]{lyons_peres_2017} implies that
\begin{align*}
\mathbb{E}\left[\lvert\overrightarrow{C_0}\rvert\right] = \sum_{x \in \Z^d}\mathbb{E}[F(0,x;\cdot)]
    =\sum_{x \in \Z^d}\mathbb{E}[F(x,0;\cdot)]=\mathbb{E}\left[\lvert\overleftarrow{C_0}\rvert\right]. 
\end{align*}
Hence, if we have out-percolation at parameters $d,k$, then at least the expected size of the in-cluster is also infinite and vice versa. A sharp-threshold result (including exponential tails of the cluster size distribution in the subcritical phase) would yield the conjecture.

Let us further remark that Conjecture~\ref{conj:in} is not known either in the continuum case. Without such a result, it seems that in-percolation is in general more difficult to deal with than out-percolation, due to increased combinatorial complexity. E.g., showing that the 1-DnG does not percolate was relatively straightforward (cf.~the proof of Proposition \ref{prop:deterministic-percolation-dng}), but proving the lack of in-percolation in the same graph (cf.~the proof of Proposition~\ref{prop-U1}) already presented more challenges. Further, in the two-dimensional Poisson case~\cite[Proof of Theorem 2]{BalisterBollobas13}, the authors showed directly that 
out-percolation occurs for $k=13$, but for in-percolation for the same $k$, the same method did not work. 
They derived that their arguments for out-percolation actually imply strong percolation and therefore also in-percolation.

We have seen that for $d$ large, out-percolation in the $k$-DnG occurs already for $k=3$, while for the $k$-BnG percolation definitely requires $k=\Omega(\sqrt d)$, and all the examples where we have verified percolation are such that $k>d$. It is an interesting open question whether strong percolation is closer to directed (out-)percolation than to bidirectional percolation in this respect.  The aforementioned arguments of~\cite{BalisterBollobas13} suggest the former. 

\section{Proofs}\label{sec:proof}
\subsection{Proofs for the $k$-DnG}\label{sec:proofs-dng}
We start by treating the cases in which the occurence of an infinite cluster that contains the origin is a deterministic event (up to $\P$-nullsets). 
\begin{proof}[Proof of Proposition \ref{prop:deterministic-percolation-dng}]
    \textit{Ad (i):} 
    If there exists a path starting from the origin and reaching an $\ell_1$-distance $n$ from the origin, then this path is unique and has at least $n$ steps. But such a path exists with probability at most $(2d)^{-n}$, which converges to zero as $n$ tends to infinity.
    \medskip
    \newline 
    \textit{Ad (ii):} 
    We will use a growth argument and show that the maximal distance to the origin is strictly increasing between generations. We initialize $G_0=\{0\}$, and construct $G_{n+1}$ from $\cup_{j=0}^nG_j$ as follows. Each vertex $v\in G_n$ independently samples $k$ neighbors uniformly at random, say $\{v_1,\ldots, v_k\}$. Then, 
    \be 
    G_{n+1}=\Big(\bigcup_{v\in G_n}\{v_1,\ldots, v_k\}\Big)\setminus \bigcup_{j=0}^n G_j. 
    \ee
    For $x\in\Z^d$, we let $\|x\|_1$ denote the $\ell_1$-distance between $x$ and the origin. We show by induction that 
    \begin{align}\label{eqn:growth-property}
         \forall n  \in \N:   \quad \max_{x \in G_n}\|x\|_1 < \max_{x \in G_{n+1}}\|x\|_1.  
    \end{align}
    This then implies that there exists an infinitely long directed path starting at the origin. 
    
    Indeed, for $n=0$ the inequality in~\eqref{eqn:growth-property} is clear. For the induction step, let $x \in G_n$ be a vertex (possibly non unique) that achieves the maximum $\ell_1$ distance. Then $x$ has $d$ neighbors $y$ with $\|x\|_1 < \|y\|_1$. If $k  \geq 2d-d+1=d+1$, then we necessarily have to choose one of these $y$ as a potential vertex for the new generation. This $y$ cannot have been in any of the previous generations by the induction hypothesis. Therefore we have $y \in G_{n+1}$ and~\eqref{eqn:growth-property} follows. 
\end{proof}

\begin{proof}[Proof of Theorem \ref{thm-CoxDurrett}]
We follow \cite[Section 2]{CoxDurrett}. Since we have already verified Proposition~\ref{prop:deterministic-percolation-dng}, without loss of generality we will assume that $k \leq d$. Let $\mathcal R_n$ denote the set of directed paths from the origin to level $n$ (that is, all vertices at $\ell_1$-distance $n$) with strictly increasing $\ell_1$-distance that stay within the first closed quadrant, in other words all vertices on the path have non-negative coordinates, and let $N_n$ be the number of open paths in $\mathcal R_n$. Then $W_n=(2/k)^nN_n$ is a martingale with respect to $\mathcal F_n$, the sigma-algebra generated by the first $n$ steps of a simple random walk $T$ with $T_0=o$ and $T_{n+1}=T_n+e_i$ where $i$ is uniformly chosen from $\{1,\dots ,d\}$ for all $n\in\N_0$. Indeed, using the fact that the expected number of monotone (i.e., with respect to $\ell_1$-distance increasing) edges  is given by the expectation of a hypergeometric random variable, i.e., 
\begin{align*}
\sum_{\ell=0}^k \ell\binom{d}{\ell}\binom{d}{k-\ell}\Big/\binom{2d}{k}=\frac{k}{2},
\end{align*}
and we can use independence to see that 
\begin{align*}
\E[N_{n+1}|\mathcal F_n]=kN_n/2.
\end{align*}

Since $\E[W_1]=1$, we can apply the martingale convergence theorem, to ensure the existence of a random variable $W$ with 
\begin{align*}
W_n\to W\qquad \text{almost surely}.
\end{align*}
Using the second-moment method it now suffices to show that $\limsup_{n\uparrow\infty}\E[W_n^2]<\infty$, since by the Paley--Zygmund inequality we then have that $\P(W>0)>0$, which implies $\theta^{\rm D}(k,d)>0$.

For this, note that 
\begin{align*}
\E[N_n^2]=\sum_{s,t\in \mathcal R_n}\P(s,t\text{ open})=\sum_{s,t\in \mathcal R_n}p^{K(s,t)}q^{L(s,t)}p^{2(n-K(s,t)-L(s,t))},
\end{align*}
where $p=k/(2d)$ denotes the probability that a given edge is open, $q=k(k-1)/(2d(2d-1))$ denotes the probability that two different edges that emerge from the same vertex are open, $K(s,t)$ is the number of joint edges in two paths $s$ and $t$, and $L(s,t)$ is the number of vertices $x$ in $s$ and $t$ such that $s$ and $t$ do not have a joint edge after $x$.  Now note that $q<p^2$, which implies that
\[ \E[N_n^2] \leq \sum_{s,t\in \mathcal R_n} p^{2n-K(s,t)}. \]
The right-hand side is precisely the same as the right-hand side of the display below \cite[Equation~(2$\cdot$7)]{CoxDurrett}. Thus, a verbatim application of the arguments of \cite[Section 2]{CoxDurrett} implies that 
\[ \lim_{n\to\infty} \E W_n^2 <\infty \]
holds if and only if 
\[ p =k/(2d)> \varrho(d), \]
where for two independent (simple, symmetric, nearest-neighbor) random walks $\widetilde S=(\widetilde S_n)_{n\in\N_0}$ and $\widetilde S'=(\widetilde S'_n)_{n\in\N_0}$ on the first quadrant of $\Z^d$ started from the origin, we define
\[ \varrho(d) = \P(\exists m\geq 0 \colon \widetilde S_m = \widetilde S'_m \text{ and } \widetilde S_{m+1} = \widetilde S'_{m+1}). \]
Let further $\tau_d=\inf \{ m \geq 0 \colon \widetilde S'_m \text{ and } \widetilde S_{m+1} = \widetilde S'_{m+1} \}$.
According to \cite[p.~155]{CoxDurrett}, we have for all $d$ that
\begin{equation}\label{firstfew} 
    \P(\tau_d=0)=d^{-1}, \qquad \P(\tau_d=1)=0, \qquad \P(\tau_d=2)=d^{-3}-d^{-4}, 
\end{equation}
and for $3 \leq l \leq d$,
\[ \P(\tau=l) \leq d^{-l} l!, \]
while for $j\geq 1$ and $l >jd$,
\begin{equation}\label{bigl} 
    \P(\tau=l) \leq d^{-1} (2\pi d)^{1/2} \Big( \e^{-1/13}/\sqrt{2\pi} \Big)^d j^{-\frac{1-d}{2}},
\end{equation}
so that
\begin{equation}\label{CDUB} 
    \begin{aligned} 
        \varrho(d)=\P(\tau_d<\infty) 
        &\leq 
        d^{-1} + d^{-3} -d^{-4} + \sum_{l=3}^d d^{-l} l! +  (2\pi d)^{1/2} \Big( \frac{\e^{-1/13}}{\sqrt{2\pi}} \Big)^d \sum_{j=1}^{\infty} j^{-\frac{1-d}{2}} 
        \\ 
        & =
        d^{-1} + d^{-3} -d^{-4} + \sum_{l=3}^d d^{-l} l! +  (2\pi d)^{1/2} \Big( \frac{\e^{-1/13}}{\sqrt{2\pi}} \Big)^d \zeta\Big( \frac{d-1}{2} \Big), 
\end{aligned}\end{equation}
where $\zeta$ denotes the Riemann zeta function. Let $U(d)$ denote the right-hand side of the inequality. To conclude that we do indeed percolate it suffices to verify that $U(d)$ is less than $k/(2d)$. In low dimensions one can compute the right-hand side numerically. 
The Table~\ref{table-CD} shows the smallest values of $k$ such that $U(d)$ is less than $k/(2d)$, for $d=4,5,6,7$.
So it remains to verify the claim for $d \geq 7$ and $k\geq 3$. For these standard but tedious calculations we refer to the proof of Lemma~\ref{lemma-largedmon}, which is given below in full detail. 
\end{proof}
\begin{table}[hb]
\begin{center}
    \begin{tabular}{|c|c|c|c|c|} \hline
       &  $d=4$ & $d=5$ & $d=6$ & $d=7$ \\ \hline
        $U(d)$, the r.h.s.\ of~\eqref{CDUB}   & 0.693093  & 0.394622  & 0.268615 & 0.199707  \\ \hline
         Smallest $k$ such that $U(d)$ is less than $k/(2d) $ & 6 & 4 & 4 & 3 \\ \hline 
    \end{tabular}
\end{center}
\caption{Values of the right-hand side of~\eqref{CDUB} and the smallest $k$ such that this right-hand side is below $k/(2d) $, for $4 \leq d \leq 7$.}\label{table-CD}
\end{table}

Before we start the technical calculations, let us briefly note that by~\eqref{firstfew} we directly have that $\varrho(d) >1/d$. Therefore, $\varrho(d)<k/(2d)$ never holds for $k=2$, whence for $k=2$ the approach of the proof of Theorem~\ref{thm-CoxDurrett} is not applicable. 
For $d\leq 3$ (where the cases $k=2,d=2$, $k=2,d=3$, and the case $k=3,d=3$ are open), the issue is that the sum $\sum_{j=1}^{\infty} j^{\frac{-1-d}{2}}$ (cf.\ \eqref{CDUB}) does not converge. 
Indeed, this technique, like many others in statistical physics, only works when we are in at least four dimensions.
Finally, for $d=4$ and $k=3$ one could hope that a combination of the proof techniques of the theorem and some explicit computations can work, but for this, one would need a sufficiently tight upper bound on $\P(\tau_4 = k)$ up to $k \approx 12$, which exceeds our available computing capacity. 

\begin{lemma}\label{lemma-largedmon}
For any $d \geq 7$ and $k \geq 3$, we have
\begin{equation}\label{largedmon} \frac{1}{d} + \frac{1}{d^3} -\frac{1}{d^4} + \sum_{l=3}^d d^{-l} l! +  (2\pi d)^{1/2} \Big( \frac{\e^{-1/13}}{\sqrt{2\pi}} \Big)^d \zeta\Big( \frac{d-1}{2} \Big) - \frac{k}{2d} < 0,  
\end{equation}
where $\zeta(\cdot)$ is the Riemann zeta function. 
\end{lemma}

\begin{proof}
    Recall that we have already seen that
    \[ 
    R(k,d)
    :=
    \frac{1}{d} + \frac{1}{d^3} -\frac{1}{d^4} + \sum_{l=3}^d d^{-l} l! +  (2\pi d)^{1/2} \Big( \frac{\e^{-1/13}}{\sqrt{2\pi}} \Big)^d \zeta\Big( \frac{d-1}{2} \Big) - \frac{k}{2d} < 0  
    \]
    holds for any $k\geq 3$ if $d=7$. 
    Our goal is to prove the statement of Lemma~\ref{lemma-largedmon}, i.e., that the same holds for any $d \geq 7$ and $k\geq 3$. We check numerically that it also holds for $8 \leq d \leq 11$ and $k=3$ (and thus also for $k\geq 4$ for the same values of $d$): indeed, we have $R(8,3)=-0.0292277$, $R(9,3)=-0.0350912$, $R(10,3)=-0.0367514$ and $R(11,3)=-0.0364418$. 

    Now,  for $d \geq 11$ we have
   \begin{equation} \label{inequ}
    \begin{aligned}
            & 
            R(k,d+1)-(d+1)^{-(d+1)}(d+1)! 
            \\ 
            &=
            \frac{1}{d+1}- \frac{k}{2(d+1)} + \frac{1}{(d+1)^3} -\frac{1}{(d+1)^4} + \sum_{l=3}^{d} \frac{ l!}{(d+1)^l} +  (2\pi (d+1))^{1/2} \Big( \tfrac{\e^{-1/13}}{\sqrt{2\pi}} \Big)^{d+1} \zeta\Big( \frac{d}{2} \Big)
            \\
            &=
            \frac{d}{d+1} \Big( \frac{1}{d}- \frac{k}{2d} \Big) + \frac{1}{(d+1)^3} -\frac{1}{(d+1)^4} + \sum_{l=3}^{d} \frac{ l!}{(d+1)^l} +  (2\pi (d+1))^{1/2} \Big( \tfrac{\e^{-1/13}}{\sqrt{2\pi}} \Big)^{d+1} \zeta\Big( \frac{d}{2} \Big)
            \\ 
           &  \leq
            \frac{d}{d+1} \Big( \frac{1}{d}- \frac{k}{2d} \Big) + \frac{d}{d+1} \bigg[ \frac{1}{d^3} -\frac{1}{d^4} + \sum_{l=3}^d d^{-l} l! +  (2\pi d)^{1/2} \Big( \tfrac{\e^{-1/13}}{\sqrt{2\pi}} \Big)^d \zeta\Big( \frac{d-1}{2} \Big) \bigg]  
            \\ & \qquad    + \frac{1}{(d+1)^3} -\frac{1}{(d+1)^4} - \frac{d}{d+1} \big( \frac{1}{d^3}-\frac{1}{d^4} \big) 
            \\ 
            &  =
            \frac{d}{d+1} R(k,d) + \Big( \frac{1}{(d+1)^3} - \frac{1}{(d+1)^4} \Big) - \frac{d}{d+1} \big( \frac{1}{d^3}-\frac{1}{d^4} \big) , 
    \end{aligned}
    \end{equation}
    where in the penultimate step we used that the Riemann zeta function is monotone decreasing on $(1,\infty)$ and that for $d \geq 11$,
    \begin{equation*} 
        \sqrt{\frac{(d+1)}{ d}} \frac{\e^{-1/13}}{\sqrt{2\pi}} \leq \sqrt{\frac{12}{11}} \frac{\e^{-1/13}}{\sqrt{2\pi}}  \approx 0.385831, 
    \end{equation*}
    while $d/(d+1)>11/12$, which is strictly larger. 

    Next, let us show that
    \begin{equation}\label{newtermbound} 
        (d+1)^{-(d+1)} (d+1)! \leq \frac{d}{d+1} \big( \frac{1}{d^3}-\frac{1}{d^4} \big) - \Big( \frac{1}{(d+1)^3} - \frac{1}{(d+1)^4} \Big) 
    \end{equation}
    holds for $d \geq 11$. Note that for such $d$, we have
    \[ 
        (d+1)^{-(d+1)} (d+1)! = \frac{(d+1)!}{12!(d+1)^{(d+1)-12}}\frac{12!}{(d+1)^{12}}\leq \frac{12!}{(d+1)^{12}} \leq \frac{12!}{12^8} \frac{1}{(d+1)^4} < \frac{1.115}{(d+1)^4}, 
    \]
    and
    \[ 
        \frac{d}{d+1} \Big( \frac{1}{d^3}-\frac{1}{d^4} \Big) - \Big( \frac{1}{(d+1)^3} - \frac{1}{(d+1)^4} \Big) = \frac{(d-1)(d+1)^3-d^4}{d^3(d+1)^4} = \frac{2d^3-2d-1}{d^3(d+1)^4} > \frac{1.99}{(d+1)^4}.
    \]
    Hence,~\eqref{newtermbound} holds for $d \geq 11$, and thus by~\eqref{inequ}, for such $d$ and for any $k\geq 3$,
    \[ 
        R(k,d+1) \leq \frac{d}{d+1} R(k,d). \numberthis\label{goodbound}
    \]
    Since $R(k,11) < 0$, the lemma now follows by induction over $d$.
\end{proof}

\begin{proof}[Proof of Theorem \ref{theorem-3322}]
    We use a randomized coupling approach. The coupling is designed to have the following two key properties. First, using the additional dimension, the edge distribution of a node $x'=(x_1,\dots,x_d)\in \Z^d$ is reproduced as the image measure of a (random) mapping from the nodes in the column $\bar x=\{(x_1,\dots,x_d,n)\in \Z^{d+1}\colon n\in \Z\}$ perpendicular to $\Z^d$ in that node and in particular thus remains iid over the nodes in $\Z^d$. Second, any connected component in $\Z^d$ is the image of a connected component in $\Z^{d+1}$.

To make this precise, let us write $x=(x_1,\dots,x_{d+1})\in \Z^{d+1}$ and $x'=\pi(x)=(x_1,\dots, x_d)\in \Z^d$ for the projection of $x$ to its first $d$ coordinates. By $\omega_{x}$ we denote the $(d+1)$-dimensional $(k+1)$-neighbor directed edge configuration of $x\in \Z^{d+1}$ and (with a slight abuse of notation) by $\omega'_{x'}=\pi(\omega_{x})$ the associated edge configuration of $x'$ in $\Z^d$. In words, we forget all edges facing up or down and only keep all remaining edges which we consider as part of the edge set of $x'$ in $\Z^d$, see Figure~\ref{fig_edgemap} for an illustration. 
\begin{figure}[h!]
\begin{tikzpicture}[scale=2]
\draw [->,ultra thick] (-2,0) -- (-2,1);
\draw [->] (-2,0) -- (-2,-1);
\draw [->] (-2,0) -- (-3,0);
\draw [->,ultra thick] (-2,0) -- (-1,0);
\draw [->,ultra thick] (-2,0) -- (-1.5,0.5);
\draw [->] (-2,0) -- (-2.5,-0.5);

\draw [->,ultra thick] (-1,0.3) to [bend left] (1,0.3);
\draw node {$\pi$};

\draw [->,ultra thick] (2,0) -- (3,0);
\draw [->] (2,0) -- (1,0);
\draw [->,ultra thick] (2,0) -- (2.5,0.5);
\draw [->] (2,0) -- (1.5,-0.5);
\end{tikzpicture}
\caption{A vertex in $\Z^3$ with $k=3$ outgoing open edges (left) and its image (right) under $\pi$ in $\Z^2$, where the open edge that is perpendicular to $\Z^2$ is removed.}\label{fig_edgemap}
\end{figure}
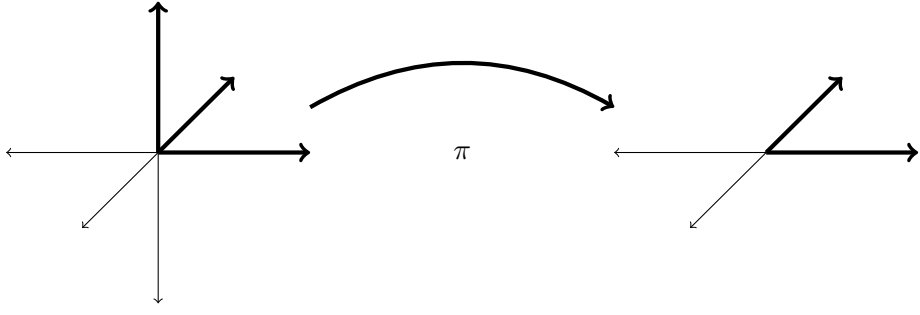

We further write $\vec{x}=(x',n)_{n\in\Z}=\pi^{-1}(x')$ for the vector of nodes perpendicular to $\Z^d$ in $x'\in\Z^d$ (in the $(d+1)$-th dimension) and $\vec\omega_{\vec x}$ for the associated vector of $(k+1)$-neighbor configurations of the nodes $(x',n)\in \Z^{d+1}$ in $\vec x$. The idea now is to define probability kernels $\nu_{x'}(\d \omega_{x'}|\vec\omega_{\vec x})$, step-by-step with respect to $x'$, in such a way that the joint kernel has two properties:
\begin{enumerate}
\item Integrating the joint kernel with respect to $\omega$ gives the distribution of the $d$-dimensional directed $k$-DnG. 
\item If the connected component of the origin in a $(d+1)$-dimensional component $\omega$ is finite, then the joint kernel only puts positive mass on finite components of the origin in $d$ dimensions as well. 
\end{enumerate}

To define this precisely, we need to consider ordered $k$-subsets of nearest neighbors of nodes $x'\in \Z^d$ with some additional coordinate, i.e., we write 
\begin{align*}
D_{x'}
=
\{((a_{1},n_1),\dots,(a_{k},n_k))\colon n_1,\dots,n_k\in \Z, \{a_1,\dots,a_k\}
\subset N^d(x')\text{ with }a_i\neq a_j \ \forall i\neq j\}
\end{align*}
where $N^d(x')\subset\Z^d$ denotes the set of nearest neighbors of $x'$ in $\Z^d$.
Then, we start with the origin $o'\in \Z^d$ and define for $\vec\omega_{\vec o}=\pi^{-1}(\omega'_{o'})$ the probability kernel $\nu_{o',0}(\cdot|\vec\omega_{\vec o})$ (at level $0$) as a probability measure on $D_{o'}$ as follows, see Figure~\ref{fig_edgemap_2} for an illustration.
\begin{enumerate}
\item[\bf Step 1] If $|\pi(\omega_{(o',0)})|\ge k$, pick $k$ neighbors $a_{o',1},\dots,a_{o',k}$ uniformly at random from the available $|\pi(\omega_{(o',0)})|$ connected neighbors in $\Z^d$. Write the outcome as $\big\{(a_{o',1},0),\dots,(a_{o',k},0)\big\}$.
\item[\bf Step 2] If $|\pi(\omega_{(o',0)})|<k$, we keep the $(k-1)$-many available connected neighbors in $\Z^d$ as\linebreak $a_{o',1},\dots, a_{o',k-1}$. Then, randomly choose a direction $\pm e_{d+1}$ (say $e_{d+1}$) perpendicular to $\Z^d$ (there must exist a directed edge facing in both perpendicular directions) and follow these steps:
\begin{enumerate}
\item If $|\pi(\omega_{(o', 1)})\setminus \pi(\omega_{(o',0)})|\ge 1$, 
pick the missing neighbor $a_{o',k}$ uniformly at random from the available connected neighbors $\pi(\omega_{(o', 1)})\setminus \pi(\omega_{(o',0)})$. Write for the total outcome $\big\{(a_{o',1},0),\dots,(a_{o',k-1},0),(a_{o',k},1)\big\}$.
\item If $|\pi(\omega_{(o', 1)})\setminus \pi(\omega_{(o',0)})|= 0$, follow the existing arrow in the same direction $e_{d+1}$ (respectively $-e_{d+1}$) as in Step~2 and repeat from (a) with $\omega_{(o',1)}$ replaced by $\omega_{(o', 2)}$ (respectively with $\omega_{(o',-1)}$ replaced by $\omega_{(o',-2)}$).
\end{enumerate}
\end{enumerate}
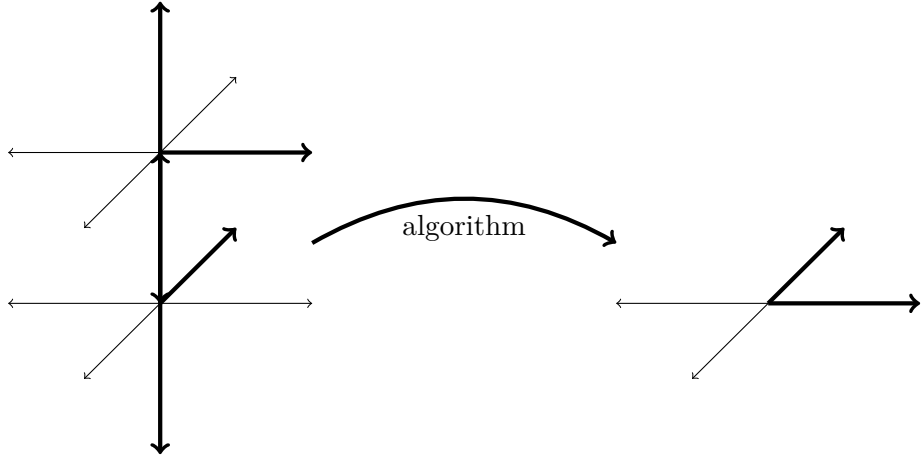
\begin{figure}[h!]
\begin{tikzpicture}[scale=2]
\draw [->,ultra thick] (-2,-0.5) -- (-2,0.5);
\draw [->,ultra thick] (-2,-0.5) -- (-2,-1.5);
\draw [->] (-2,-0.5) -- (-3,-0.5);
\draw [->] (-2,-0.5) -- (-1,-0.5);
\draw [->,ultra thick] (-2,-0.5) -- (-1.5,0);
\draw [->] (-2,-0.5) -- (-2.5,-1);

\draw [->,ultra thick] (-2,0.5) -- (-2,1.5);
\draw [->, ultra thick] (-2,0.5) -- (-2,-0.5);
\draw [->] (-2,0.5) -- (-3,0.5);
\draw [->,ultra thick] (-2,0.5) -- (-1,0.5);
\draw [->] (-2,0.5) -- (-1.5,1);
\draw [->] (-2,0.5) -- (-2.5,0);

\draw [->,ultra thick] (-1,-0.1) [bend left] to (1,-0.1);
\draw node {algorithm};

\draw [->,ultra thick] (2,-0.5) -- (3,-0.5);
\draw [->] (2,-0.5) -- (1,-0.5);
\draw [->,ultra thick] (2,-0.5) -- (2.5,00);
\draw [->] (2,-0.5) -- (1.5,-1);
\end{tikzpicture}
\caption{Two vertices $(o',0)$ and $(o',1)$ in $\Z^3$  with $k+1=3$ outgoing open edges (left). According to Step~2a in the algorithm (Step~1 is not applicable), we first keep the open edges of $(o',0)$, which are inside $\Z^2$, and then (according to Step~2b, at random) we move upwards to $(o',1)$, where we find the second (missing) edge such that $o'\in\Z^2$ has now degree $k=2$.}\label{fig_edgemap_2}
\end{figure}
Now, almost surely, the construction must end (since the probability for an infinite repetition of Step~2 is zero) and, due to the complete symmetry in the construction, $\nu_{o',0}$ reproduces the uniform $k$-nearest-neighbor distribution for the origin. That is, for 
$$E^{k,d}_{x'}=\{\{a_1,\dots, a_k\}\subset N^d(x')\}$$ 
and $F\colon D_{x'}\to E^{k,d}_{x'}$, $\big\{(a_1,n_1),\dots,(a_k,n_k)\big\}\mapsto \{a_1,\dots, a_k\}$ and $U^{k,d}_x$ the uniform distribution on $E^{k,d}_x$, we have that, for all $\omega'_o\in E^{k,d}_o$,
\begin{align*}
\int \Big(\bigotimes_{n\in \Z}U^{k+1,d+1}_{(o',n)}\Big)(\d \vec\omega_{\vec o})\nu_{o',0}(F^{-1}(\omega'_{o'})|\vec\omega_{\vec o})=U^{k,d}_{o'}(\omega'_{o'}).
\end{align*}

Let us describe the following steps in words. Under $\nu_{o',0}(\cdot|\vec\omega_{\vec o})$ we are equipped with $k$ neighbors of the origin in $\Z^d$ that all also carry the information on which level in the additional $(d+1)$-th coordinate they are discovered. Starting with $\hat a_{o',1}=(a_{o',1},n_1)$, the first discovered neighbor, we can sample its neighbors, using the same algorithm, based on the information provided by $\bar a_{o',1}$, now started at level $n_1$ (and not at level $0$). The same can be done for all other neighbors $\hat a_{o,i}=(a_{o,i},n_i)$. In this fashion, we can step-by-step explore the connected component of the origin in $\Z^d$, where in any step, whenever we discover a new vertex $y\in \Z^d$, we only use information from the associated vector $\vec y$, i.e., 
\begin{align*}
\nu(\d \omega'|\omega)=&\int \nu_{(o',0)}(\d(\hat a_{o',1},\dots,\hat a_{o',k})|\vec \omega_{\vec o})\times\\
&\int\nu_{(a_{o',1},n_1)}(\d(\hat a_{a_{o',1},1},\dots,\hat a_{a_{o',1},k})|\vec \omega_{\vec a_{o',1}})\dots \int\nu_{(a_{o',k},n_k)}(\d(\hat a_{a_{o',k},1},\dots,\hat a_{a_{o',k},k})|\vec \omega_{\vec a_{o,k}})\times\\
&\cdots \1\{F(\hat a_{o',1},\dots,\hat a_{o',k},\hat a_{a_{o',1},1},\dots,\hat a_{a_{o',1},k},\dots)=\d \omega'\}.
\end{align*}

Now, the key point is the following. Imagine that the algorithm presents us a configuration of open edges in $\Z^d$ such that the connected component of the origin reaches the boundary of a centered box with side-length $n$. Then, necessarily, any configuration in $\Z^{d+1}$ that is used as an input for the algorithm must also have the property that the origin is connected to the boundary of a box of side-length $n$, now of course in  $\Z^{d+1}$. Hence, denoting by $A^{d}_n$ the event that the origin is connected to the boundary of the centered box of side-length $n$ in $\Z^d$ and by $\P_{k,d}$ the distribution of the $k$-DnG in dimension $d$, we have that 
\begin{align*}
\P_{k,d}(A^d_n)=\E_{k+1,d+1}[\nu(A^d_n|\cdot)]\le \P_{k+1,d+1}(A^{d+1}_n), 
\end{align*}
which gives the result.
\end{proof}

\subsection{Proofs for the $k$-UnG}\label{sec:kungproof}
\begin{proof}[Proof of Theorem \ref{theorem-U22}]
    We prove the positive percolation probability in the 2-UnG via a dual approach. That is, let $\Z'^2:=\{x+(1/2,1/2): x\in \Z^2\}$ be the two-dimensional dual lattice of $\Z^2$. An edge $e'\in\Z'^2$ has exactly one edge $e\in\Z^2$ that crosses $e'$, and we say the edge $e'$ is open if and only if the edge $e$ is open. For the $2$-UnG model in two dimensions it thus follows that for any such edge $e'\in\Z'^2$, 
    \be 
    \P(e' \text{ is closed})=1-\P(e\text{ is open})=1-3/4=1/4.
    \ee 
    Moreover, the negative correlations between edges that are present in the $2$-UnG model also appear in its dual. Namely, whenever for edges $e'_1,e'_2\in\Z'^2$ their unique associated edges $e_1,e_2\in\Z^2$ are neighbors, then
    \be 
    \P(e'_1,e'_2\text{ closed})=1/24\leq 1/16=\P(e'_1\text{ closed})^2.
    \ee 
    Otherwise, the status of the two dual edges is independent.
    This can be extended, so that for any $k\geq 2$ and any self-avoiding path $(e'_1,\dots, e'_k)$ of nearest-neighbor vertices in $\Z'^2$ with $\P(e'_1,\dots, e'_k\text{ closed})>0$, 
    \be 
    \P(e'_1,\dots, e'_k\text{ closed})= \prod_{i=0}^{k-1}\P(e'_{i+1}\text{ closed}|e'_1,\dots,e'_i\text{ closed})\le 4^{-k}.
    \ee 
    As such, we can use the negative correlation to bound the probability that there exists a closed circuit in the dual that surrounds the origin. For any $n\geq 4$, the number of such circuits of length $n$ is at most $nc_{n-1}(2)$, where we recall that $c_n(2)$ denotes the number of self-avoiding paths of length $n$ in $\Z^2$. As $c_n(2)=(c(2)+o(1))^{n-1}$ by the definition of $c_n$ and $c$ (cf.~Section~\ref{sec:B}), we thus arrive at 
    \be \ba\label{eq:finnorpercprob}
    1-\theta^{\rm U}(2,2)&=\P(\exists \text{ closed circuit }\gamma \text{ in }\Z'^2\text{ that surrounds the origin})\\
    &\leq \sum_{n\ge4} \sum_{\text{ circuits }\gamma\colon|\gamma|=n}\P(\gamma\text{ closed in }\Z'^2)\leq \sum_{n\ge4} n(c(2)+o(1))^{n-1} 4^{-n}, 
    \ea \ee 
    which is finite since $c(2)\leq 3$. We now apply an argument from~\cite{JanSwart} (see also the references therein) to show that this probability is in fact smaller than one, which yields $\theta^{\rm U}(2,2)>0$. Let $D_m=\{0,\ldots, m\}^2$ and say that a vertex $i\in \Z^2$ is `wet' when there exists a $j\in D_m$ such that $i\to j$, that is, when there exists a path of open edges from $i$ to $j$ in the primary lattice. Suppose the component of the origin is finite almost surely. Then, the number of wet vertices is finite almost surely as well, and there exists a circuit of closed edges in the dual that surrounds all the wet sites. By the same argument as in~\eqref{eq:finnorpercprob}, 
    \be 
    \P(D_m \not\to \infty )\leq \sum_{n\ge4m} \sum_{\text{ circuits }\gamma\colon |\gamma|=n}\P(\gamma\text{ closed in }\Z'^2)\leq \sum_{n\ge4m} n(c(2)+o(1))4^{-n},
    \ee 
    where we can start the outer sum from $4m$ since a circuit that surrounds all wet vertices must surround $D_m$ and hence have length $4m$ at least. 
    As a result, since the sum is finite for each $m\in\N$, choosing $m$ large enough yields that $\P(D_m\to \infty)>0$ and hence $\P(i\to \infty)>0$ for some $i\in D_m$. By translation invariance, it then follows that $\theta^{\rm U}(2,2)>0$, as desired.

    To prove that $\theta^{\rm U}(3,3)>0$, it suffices to show that with positive probability $o$ is connected to infinity via any two-dimensional plane $S$ of $\Z^3$ including $o$. This can be done via the dual approach employed in the first part of the proof. 
    Namely, each edge $e$ of the $3$-UnG is open with probability
    \be 
    \P(e \text{ open})=1-\P(e\text{ closed})=1-\bigg(\binom{5}{3}\Big/ \binom{6}{3}\bigg)^2=\frac34,
    \ee
    which applies in particular to all edges contained in the restriction of the $3$-UnG to $S$. Thus, for any edge $e'$ of the dual lattice corresponding to $S$, we have $\P(e'\text{ closed})=1/4$. Then, for any two dual edges for which their associated non-duals are not nearest neighbors, their status is independent, and otherwise,
    \[ 
        \P(e'_1,e'_2\text{ closed}) = \binom{5}{3}\binom{5}{3}\binom{4}{3}\Big/\binom{6}{3}^3 = \frac{1}{20}\leq \frac{1}{16} = \P(e'_1\text{ closed})^2,
    \]
whence the proof can be finished analogously to the one for the $2$-UnG. 
\end{proof}

\begin{proof}[Proof of Proposition \ref{prop-U1}]
    Assume for a contradiction that we do have $\theta^{\rm U}(1,d)>0$ for some $d \in \N$, and let $A$ then denote the event of positive probability that there exists an infinite path starting from $o$ consisting of edges of the $1$-UnG. 
Let $o=X_0,X_1,X_2,\ldots$ be any infinite path of the $1$-UnG starting from $o$.
Now, if $X_0,X_1,X_2$ is the sequence of vertices in an infinite path (on the event $A$), we put
\[ K=\inf \{ k \geq 0 \colon (X_k,X_{k+1}) \text{ is not an edge of the $1$-DnG} \}. \]
We know from the proof of Part $(i)$ in Proposition~\ref{prop:deterministic-percolation-dng} that $K<\infty$ almost surely. Next, let us define
\[ L=\inf \{ k > K \colon 
(X_{k},X_{k+1})\text{ is an edge of the 1-DnG}\}. \]
Now we claim that $L=\infty$ on the event $\{ K<\infty\}$ (where we put $K=\infty$ on the event $A^c$). Indeed, for all $k=K,K+1,\ldots,L-1$, $(X_{k+1},X_k)$ is an edge of the $1$-DnG because it follows from the definition of $K$ that $(X_{k},X_{k+1})$ is not an edge of the 1-DnG, but $(X_k,X_{k+1})$ is an edge of the 1-UnG, and this is only possible if $(X_{k+1},X_k)$ is an edge of the $1$-DnG. But now, if $L<\infty$, then $(X_{L},X_{L-1})$ is an edge of the $1$-DnG, so that $(X_{L},X_{L+1})$ cannot be an edge of the $1$-DnG since there is only one edge going out of $X_L$. This implies the claim.

Let us denote the $\ell_1$-sphere of radius $n\in\N_0$ by $S^{(1)}_n$ (with $S^{(1)}_0=\{o\}$). Note that we have that $\P(\kappa<\infty\,|\,A)=1$, where
\[ \kappa = \inf \big\{ k \geq 0 \colon o \not\rightsquigarrow S^{(1)}_k \text{ in the $1$-DnG} \big\}, \]
represents the event that there exists no path along open (directed) edges from the origin to some vertex in $S^{(1)}_k$. 
Since $o \in S_0^{(1)}$ and $o$ is always connected to one of its nearest neighbors by a directed edge, we have $\P(\kappa \geq 2)=1$.
Given that $\kappa$ is an $\N_0 \cup \{\infty\}$-valued random variable, we can find $k_0\in\N \setminus \{ 1 \}$ and $\varepsilon>0$ such that $\P(\kappa =  k_0 |A) \geq \varepsilon$.  Now, by our previous observation that $L=\infty$ on $\{ K<\infty\}$, since from $o$ to $S^{(1)}_{\kappa-1}$ there always exists a directed path in the $1$-DnG, 
on the event $\{ \kappa = k_0 \} \cap A$, there must exist a self-avoiding directed path of length $n$ in the $1$-DnG ending at some vertex of $S^{(1)}_{k_0-1}$  for all $n\in\N_0$. 

Thus, writing $c_{n,k_0-1}(d)$ for the number of self-avoiding paths in the $d$-dimensional lattice of length $n$ starting from (or ending at) $S^{(1)}_{k_0-1}$, it is clear that
\[ \lim_{n\to\infty} c_{n,k_0-1}(d)^{1/n} = \lim_{n\to\infty} c_{n,1}(d)^{1/n}=c(d) \leq 2d-1, \]
where $c(d)$ is the connective constant of $\Z^d$. Now, any self-avoiding path of length $n$ ending at $S^{(1)}_{k_0-1}$ is open with probability $(2d)^{-n}$. Hence,
\[ \P(\exists\text{ self-avoiding path of length $n$ starting from $S^{(1)}_{k_0-1}$ included in $1$-DnG}) \leq \Big( \frac{2d-1}{2d} \Big)^n, \]
which decays exponentially fast in $n$. This contradicts the assertion that the left-hand side is bounded from below by $\P(A)\varepsilon>0$ independently of $n$. Therefore, $\P(A)=0$, as desired.
\end{proof}

\subsection{Proofs for the $k$-BnG}
\begin{proof}[Proof of Lemma \ref{lemma-Bnoperc}]
Note first that the system is negatively correlated in the sense that for all directed edges $\vec{\ell}_1=(x_1,x_2)\neq\vec{\ell}_2=(y_1,y_2)$ we have that
$$
\mathbb P(\vec{\ell}_1\text{ open and }\vec{\ell}_2\text{ open})\le \mathbb P(\vec{\ell}_1\text{ open})^2.
$$
Indeed, if $x_1\neq y_1$, then the inequality is an equality by independence. However, if $x_1=y_1$ and $x_2\neq y_2$, then 
$$
\mathbb P(\vec{\ell}_1\text{ open}\,|\,\vec{\ell}_2\text{ open})=(k-1)/(2d-1)\le k/(2d)=\mathbb P(\vec{\ell}_1\text{ open}).
$$
Note that the negative correlation carries over also to the bidirectional case since for $\ell_1,\ell_2$ undirected with $x_1=y_1$ and $x_2\neq y_2$, then 
$$
\mathbb P(\ell_1\text{ open}\,|\,\ell_2\text{ open})=k(k-1)/(2d(2d-1))\le k^2/(2d)^2=\mathbb P(\ell_1\text{ open}).
$$
Let us denote by $\ell_n$ the straight line starting at the origin up to the node $(ne_1)$. Then, using the first moment method and negative correlations, we have that for all $n\ge 1$,
\begin{align*}
\P(o\leftrightsquigarrow \infty \text{ in the $k$-BnG on }\Z^d)&\le \P(\exists\text{ self-avoiding open path starting at }o\text{ of length }n)\\
&\le c_n(d)^n\P(\ell_n\text{ is open}).
\end{align*}
Finally, for a suitable constant $C>0$, we have that $\P(\ell_n\text{ is open})\le C p^n$ where $p=k(k-1)/(2d(2d-1))$ is the probability that the origin samples two fixed edges to be open. This leads to the criterion for absence of percolation $c(d)p<1$.
\end{proof}

\begin{proof}[Proof of Proposition~\ref{prop-Bperc}]
As in the proof of the assertion $\theta^{\rm U}(3,3)>0$ of Theorem~\ref{theorem-U22}, it suffices to verify that there is an infinite path starting from $o$ in the restriction of the $k$-BnG to a fixed two-dimensional plane $S$ including $o$ with positive probability. Using the duality approach of the proof of Theorem~\ref{theorem-U22}, dual edges of the two-dimensional lattice of $S$ are open with probability  $(k/(2d))^2$ (which is the probability that the corresponding edge of the original graph is open in the bidirectional sense).
Further, the indicators of two dual edges being open, i.e., the corresponding edges of the $k$-BnG restricted to $S$ being open, are non-positively correlated. Indeed, we saw already in the proof of Lemma~\ref{lemma-Bnoperc} that the indicators of two edges of the $k$-BnG being open are either independent or strictly negatively correlated. Now, if the indicators of two events are independent or strictly negatively correlated, then so are the indicators of the complements of the two events. We see from the proof of Theorem~\ref{theorem-U22} that given these non-positive correlations, it suffices to choose $k$ and $d$ in such a way that a dual edge is closed with probability at most $1/c(2)$. This holds whenever 
\be  
\big(k/(2d)\big)^2 > 1-1/c(2),\qquad \text{or, equivalently,}\qquad k>d\sqrt{4(1-1/c(2))}, 
\ee
as asserted. 
\end{proof}

\begin{proof}[Proof of Proposition \ref{prop-1indep}]
Note that the $k$-BnG model (just as the $k$-UnG model) in any dimension $d$ is an 1-dependent  bond percolation model where each edge is open with probability $p=k^2/(4d^2)$. Let us write 
\[ p_{\sup}(\Z^d)= \sup \{ p \in (0,1) \colon \text{some 1-dependent bond percolation model does not percolate} \}. \]
According to \cite{improvedLB}, we have 
 $\lim_{d\to\infty} p_{\sup}(\Z^d) \leq 0.5847$, which implies the statement.
\end{proof}

Let us remark that another result by \cite{improvedLB} is that $p_{\sup}(\Z^2) \leq 0.8457$. This together with the assertion that $p_{\sup}(\Z^d) \leq p_{\sup}(\Z^{d-1})$ yields that for $\alpha>2\sqrt{0.8457}$ we have $\theta^{\rm B}(\lfloor \alpha d \rfloor, d)>0$. However, this statement is weaker than Proposition~\ref{prop-Bperc} given that $2\sqrt{0.8457} \approx 1.8392 > \sqrt{4(1-1/2.679192495)}$. 

\section*{Acknowledgements}
The authors thank D.~Mitsche, G.~Pete, and B.~Ráth for interesting discussions and comments. BJ and JK acknowledge the financial support of the Leibniz Association within the Leibniz Junior Research Group on \textit{Probabilistic Methods for Dynamic Communication Networks} as part of the Leibniz Competition. BL was supported by the grant GrHyDy ANR-20-CE40-0002. AT  was partially supported by the ERC Consolidator Grant 772466 ``NOISE''.

\bibliographystyle{amsalpha}
\bibliography{bib}
\end{document}